\documentclass [a4paper,twoside,11pt]{article}

\usepackage[latin1]{inputenc}
\usepackage{amsmath,amsfonts,amssymb}
\usepackage{vmargin,graphicx,theorem}
\usepackage[english]{babel}
\usepackage{enumerate}

\setpapersize[portrait]{A4}
\setmarginsrb{2.2cm}{2cm}{2.2cm}{1cm}{0cm}{0.1cm}{0.5cm}{2cm}

\newcommand{\disp}{\displaystyle}

\newcommand{\dN}{\ensuremath{\mathbb{N}}}

\newcommand{\dR}{\ensuremath{\mathbb{R}}}

\newtheorem{ethm}{Theorem}[section]

\newtheorem{ecor}[ethm]{Corollary}

\newtheorem{elem}[ethm]{Lemma}

\newtheorem{exe}[ethm]{Example}

\newtheorem{erem}[ethm]{Remark}

\newcommand{\proofend}{~$\rhd$}
\newcommand{\proofbegin}{~$\lhd$}

\newenvironment{eproof}
               {\noindent {\emph{\textbf{Proof}}}\\\proofbegin~}
               {\proofend\\}

\newcommand{\p}[4]{{#3}\!\left#1{#4}\right#2}

\newcommand{\ABS}[1]{\ensuremath{{\left| #1 \right|}}} 
\newcommand{\PAR}[1]{\ensuremath{{\left(#1\right)}}} 
\newcommand{\BRA}[1]{\ensuremath{{\left\{#1\right\}}}} 
\newcommand{\NRM}[1]{\ensuremath{{\left\Vert #1\right\Vert}}} 

\renewcommand{\phi}{\varphi}

\renewcommand{\geq}{\geqslant}

\newcommand{\varf}[1]{\mathbf{Var}_{#1}}
\newcommand{\entf}[1]{\mathbf{Ent}_{#1}}

\newcommand{\ent}[2]{\p(){\entf{#1}}{#2}}

\newcommand{\var}[2]{\p(){\varf{#1}}{#2}}

\def\disp{\displaystyle}

\newcommand{\Ga}{\boldsymbol{\Gamma}}
\newcommand{\gu}{\boldsymbol{\Ga}}
\newcommand{\gd}{{\bf \gu {\!\!_2}}}

\newcommand{\al}{\alpha}

\newcommand{\ep}{\epsilon}

\newcommand{\la}{\lambda}

\newcommand{\e}{\varepsilon}

\newcommand{\affil}[1]{{\small\sl #1}}
\newcommand{\email}[1]{{\small E-mail: {\textsf {#1}}}}
\newcommand{\http}[1]{{\small Internet: {\textsf {#1}}}}

\begin{document}

\title{\sl Logarithmic Sobolev inequality for log-concave measure from Pr\'ekopa-Leindler inequality}
\author{
  Ivan Gentil\\
\affil{Ceremade (UMR CNRS no. 7534), Universit\'e Paris IX-Dauphine,}\\
\affil{Place de Lattre de Tassigny, 75775 Paris C\'edex~16, France}\\
\email{gentil@ceremade.dauphine.fr}\\
\http{http://www.ceremade.dauphine.fr/
\raisebox{-4pt}{$\!\!\widetilde{\phantom{x}}$}gentil/}\\
}
\date{\today}\maketitle\thispagestyle{empty}

\begin{abstract}
We develop in this paper an amelioration of the method given by S.
Bobkov and M. Ledoux in~\cite{bobkov-ledoux1}. We prove by
Pr\'ekopa-Leindler Theorem an optimal modified logarithmic Sobolev
inequality adapted for all log-concave measure on $\dR^n$. This
inequality implies results proved by Bobkov and Ledoux, 
the  Euclidean Logarithmic Sobolev inequality  generalized in the
last years and it also implies some convex logarithmic Sobolev
inequalities for large entropy.

\begin{center}
{\bf Résumé}
\end{center}

Dans cet article nous proposons une am\'elioration de la m\'ethode 
d\'evelopp\'ee par S. Bobkov et M. Ledoux dans~\cite{bobkov-ledoux1}. Nous prouvons par le 
th\'eor\`eme 
de Pr\'ekopa-Leindler une inégalit\'e de Sobolev logarihmique, optimale et adapt\'ee \`a toutes les
mesures 
log-concaves sur $\dR^n$. Cette in\'egalité implique les r\'esultats de Bobkov et Ledoux, les 
in\'egalit\'es de  Sobolev logarithlmique de type Euclidien généralis\'ees ces derni\`eres ann\'ees
et enfin cetaines inégalit\'es 
 de Sobolev logarithmiques de type convexe pour les grandes entropies.  
\end{abstract}

\section{Introduction}
\label{sec-int}

Pr\'ekopa-Leindler is the functional form of Brunn-Minkowski
inequality. Let $a,b>0$, $a+b=0$, 
and $u$, $v$, $w$ three non
negative measurable functions on $\dR^n$. Assume that, for any
$x,y \in \dR ^n$, we have
\begin{equation*}
 u(x)^a v(y)^b\leq w(a x + b y),
\end{equation*}
then
\begin{equation}
\label{ul-6.5} \left(\int u(x) dx\right)^a\left(\int v(x)
dx\right)^b\leq \int w(x) dx .
\end{equation}
If you applied inequality~\eqref{ul-6.5} to characteristic
functions of bounded measurable sets $A$ and $B$ in $\dR^n$, it
yields the multiplicative form of the Brunn-Minkowski inequality
$$
vol(A)^a vol(B)^b\leq vol(aA+bB),
$$
where $aA+bB=\BRA{ax_A+bx_B,\,\, x_A\in A,x_B\in B}$. One can see
for example two interesting reviews on this topic \cite{gupta,maurey}.

Bobkov and Ledoux in \cite{bobkov-ledoux1} use Pr\'ekopa-Leindler
Theorem to prove some functional inequalities like Brascamp-Lieb,
Logarithmic Sobolev and Transportation inequalities.

\medskip
More precisely,  let $\phi$ be a $\mathcal C^1$ strictly convex
function on $\dR^n$ and let
$$
d\mu_\phi(x)=e^{-\phi(x)}dx,
$$
the probability measure on $\dR^n$ (assume that $\int
e^{-\phi(x)}dx=1$). Bobkov-Ledoux prove in particular the following two results:
\begin{itemize}
\item (Proposition~2.1) Brascamp-Lieb inequality: assume that $\phi$ is a $\mathcal C^2$ function
then for all smooth enough $g$,
\begin{equation}
\label{eq-br} \var{\mu_\phi}{g}:=\int \PAR{g-\int
gd\mu_\phi}^2d\mu_\phi\leq\int \nabla
g\cdot\text{Hess}(\phi)^{-1}\nabla gd\mu_\phi,
\end{equation}
where $\text{Hess}(\phi)^{-1}$ is the inverse of the Hessian of
$\phi$.
\item (Proposition~3.1) Assume that for some $c>0$ and $p\geq2$,
 all $t,s>0$ with $t+s=1$, and for all $x,y\in\dR^n$, $\phi$
satisfies
\begin{equation}
\label{e-in}
t\phi(x)+s\phi(y)-\phi(tx+sy)\geq\frac{c}{p}(s+o(s))\NRM{x-y}^p,
\end{equation}
where $\NRM{\cdot}$ is the Euclidean norm in $\dR^n$. Then for all
smooth enough function~$g$,
\begin{equation}
\label{e-bl1} \ent{\mu_\phi}{e^g}:=\int e^g\log \frac{e^g}{\int
e^gd\mu_\phi}d\mu_\phi\leq c\int \NRM{\nabla g}^qe^gd\mu_\phi,
\end{equation}
where $1/p+1/q=1$. They give the example of the function $\phi(x)=\NRM{x}^p+Z_\phi$ ($Z_\phi$ is 
a normalization constant) which
satisfies inequality~\eqref{e-in} for some constant $c$.
\end{itemize}

In this article, we prove also  with Pr\'ekopa-Leindler Theorem, some
optimal logarithmic Sobolev inequality for log-concave measure
without conditions like inequality~\eqref{e-in}. We obtain, for
all smooth enough  function~$g$ on $\dR^n$,
\begin{equation}
\label{e-theo0} \ent{\mu_\phi}{e^g}\leq \int \BRA{x\cdot\nabla
g(x)-\phi^*(\nabla\phi(x))+\phi^*\PAR{\nabla\phi(x)-\nabla
g(x)}}e^{g(x)}d\mu_\phi(x),
\end{equation}
where $\phi^*$ is the Frenchel-Legendre transform of $\phi$,
$\phi^*(x):=\sup_{z\in\dR^n}\BRA{x\cdot z-\phi(z)}$.

\medskip

The $\gd$-criterion of Bakry-Emery  implies that if
$\text{Hess}\PAR{\phi}\geq\la Id$ in the sense of symmetric matrix
 with $\la>0$, then the probability measure $\mu_\phi$ satisfies classical
 logarithmic Sobolev inequality, for all smooth function $g$,
\begin{equation}
\label{eq-g}
  \ent{\mu_\phi}{e^g}\leq \frac{1}{2\la}\int \NRM{\nabla g}^2e^gd\mu_\phi.
\end{equation}
This inequality is proved by Gross in~\cite{gross}, one can see
also~\cite{logsob} for a review about this inequality and the
related fields. Inequality~\eqref{e-theo0} is then a
generalization of the classical logarithmic Sobolev inequality of
Gross, adapted for all log-concave measure on $\dR^n$ which does'nt
satisfies $\gd$-criterion. We get an optimal
 modified logarithmic Sobolev
inequality for log-concave measures. 

\bigskip

The next section is divided into two subsections. In the first one
we give the main theorem of this paper:
inequality~\eqref{e-theo0}.  In the second subsection we explain
how the theorem implies results of~\cite{bobkov-ledoux1}. In
particular one find again Brascamp-Lieb inequality~\eqref{eq-br}
or modified logarithmic Sobolev inequality for some function
$\phi$, inequality~\eqref{e-bl1}. In section~\ref{sec-e} we prove
that inequality~\eqref{e-theo0} is equivalent to the Euclidean
logarithmic Sobolev inequality. In particular it gives a short
proof of the generalization given in
\cite{del-dol,gentil03,ghoussoub}. In section~\ref{sec3} we give a
convex inequality for large entropy. In particular we obtain a
$n$-dimensional version for large entropy of inequalities prove in
\cite{ ge-gu-mi,ge-gu-mi2}.
\section{Logarithmic Sobolev inequality}
\label{sec-theo}

\subsection{The main theorem}
\label{s-re}

\begin{ethm}
\label{the-theo}
Let $\phi$ be a $\mathcal C^1$ strictly convex function on $\dR^n$, such that
\begin{equation}
\label{eq-propphi}
\lim_{\ABS{x}\rightarrow\infty}\frac{\phi(x)}{\NRM{x}}=\infty.
\end{equation}
We note the probability measure
$$\mu_\phi(dx)=e^{-\phi(x)}dx,$$
where $dx$ is the Lebesgue measure on $\dR^n$, assume that $\int e^{-\phi(x)}dx=1$.
\medskip

Then for all function $g$ on $\dR^n$,  smooth enough such that integrals used exits we have
\begin{equation}
\label{eq-theo} \ent{\mu_\phi}{e^g}\leq \int \BRA{x\cdot\nabla
g(x)-\phi^*(\nabla\phi(x))+\phi^*\PAR{\nabla\phi(x)-\nabla
g(x)}}e^{g(x)}d\mu_\phi(x).
\end{equation}
\end{ethm}

\begin{elem}
\label{lem-technique}
Let $g$ be a $\mathcal C^\infty$ function with a compact support on $\dR^n$. Let $s,t\geq 0$ with $t+s=1$ and
we note for $z\in\dR^n$,
$$
g_s(z)=\sup_{z=tx+sy}\PAR{g(x)-\PAR{t\phi(x)+s\phi(y)-\phi(tx+sy)}}.
$$
Then we get
\begin{multline*}
g_s(z)=g(z)+s\BRA{z\cdot\nabla
g(z)-\phi^*\PAR{\nabla\phi(z)}+\phi^*\PAR{\nabla\phi(x)-\nabla
g(x)}}\\+O\PAR{((z-y_0)\cdot\nabla (g+\phi)(z)+\NRM{z-y_0}^2)s^2},
\end{multline*}
where $y_0\in\dR^n$.
\end{elem}

\begin{eproof}
Let $s\in]0,1/2[$ and note $x=z/t-(s/t)y$, hence
$$
g_s(z)=\phi(z)+\sup_{y\in\dR^n}\PAR{g\PAR{\frac{z}{t}-\frac{s}{t}y}-t\phi\PAR{\frac{z}{t}-\frac{s}{t}y}-s\phi(y)}.
$$
Due to the fact that $g$ has a compact support and by the property~\eqref{eq-propphi} there exists $y_s\in\dR^n$ such that
$$
\sup_{y\in\dR^n}\PAR{g\PAR{\frac{z}{t}-\frac{s}{t}y}-t\phi\PAR{\frac{z}{t}-\frac{s}{t}y}-s\phi(y)}=
g\PAR{\frac{z}{t}-\frac{s}{t}y_s}-t\phi\PAR{\frac{z}{t}-\frac{s}{t}y_s}-s\phi(y_s).
$$
Moreover $y_s$ satisfies
\begin{equation}
\label{eq-ys}
\nabla g\PAR{\frac{z}{t}-\frac{s}{t}y_s}-t\nabla\phi\PAR{\frac{z}{t}-\frac{s}{t}y_s}+t\nabla \phi(y_s)=0.
\end{equation}
The function $\phi$ is a strictly convex function then there is a unique
solution $y_0$ of the equation
\begin{equation}
\label{eq-der} \nabla\phi(y_0)={\nabla\phi(z)-\nabla g(z)},\quad
y_0=\PAR{\nabla\phi}^{-1}\PAR{\nabla\phi(z)-\nabla g(z)}.
\end{equation}

We prove now that $\lim_{s\rightarrow 0}y_s=y_0$.

\medskip

First we prove that there exists $R\geq0$ such that $\forall
s\in[0,1/2]$, $\NRM{y_s}\leq R$. Indeed, if the function $y_s$ is
not bounded one can found $(s_k)_{k\in\dN}$ such that
$s_k\rightarrow0$ and $\NRM{y_{s_k}}\rightarrow\infty$. By
property~\eqref{eq-propphi}
$\lim_{\NRM{x}\rightarrow\infty}\phi(x)=\infty$ then since $g$ is
bounded we obtain $s_ky_{s_k}=O(1)$. Due to to the strictly
convexity of $\phi$, the last assertion is in contradiction with
equation~\eqref{eq-ys}.

Let $\hat{y}$ a value of adherence at $s=0$ of the function $y_s$
then $\hat{y}$ satisfies equation~\eqref{eq-der}. By unicity of
the solution of~\eqref{eq-der} we get  $\hat{y}=y_0$. Then we have
proved that $\lim_{s\rightarrow 0}y_s=y_0$.
\medskip

By Taylor formula and the continuity of $y_s$ at $s=0$ we get
$$
\phi\PAR{\frac{z}{t}-\frac{s}{t}y_s}=\phi(z)+s(z-y_0)\cdot\nabla\phi(z)
+O\PAR{((z-y_0)\cdot\nabla\phi(z)+\NRM{z-y_0}^2)s^2},
$$
and
$$
g\PAR{\frac{z}{t}-\frac{s}{t}y_s}=g(z)+s(z-y_0)\cdot\nabla g(z)
+O\PAR{((z-y_0)\cdot\nabla g(z)+\NRM{z-y_0}^2)s^2}.
$$
Then
\begin{multline*}
g_s(z)=g(z)+s\BRA{\phi(z)-\phi(y_0)+(z-y_0)\cdot(\nabla
g(z)-\nabla\phi(z))}\\+ O\PAR{((z-y_0)\cdot\nabla
(g+\phi)(z)+\NRM{z-y_0}^2)s^2}.
\end{multline*}
Using equation~\eqref{eq-der} and the expression of the
Frenchel-Legendre transformation for a strictly convex function
$$
\phi^*(x)=x\cdot\PAR{\nabla\phi}^{-1}(x)-\phi\PAR{\PAR{\nabla\phi}^{-1}(x)},
$$
and
$$
\phi^*(\nabla \phi (z))=\nabla\phi(z)\cdot z-\phi(z),
$$
we get the result.
\end{eproof}

{\noindent {\emph{\textbf{Proof of
Theorem~\ref{the-theo}}}\\\proofbegin~} The proof is based on the
proof of Theorem~3.2 of \cite{bobkov-ledoux1}. First we prove
inequality~\eqref{eq-theo} for all function $g$, $\mathcal
C^\infty$ with a compact support on $\dR^n$.

 Let $t,s\geq 0$ with $t+s=1$ and we note for $z\in\dR^n$,
$$
g_t(z)=\sup_{z=tx+sy}\PAR{g(x)-\PAR{t\phi(x)+s\phi(y)-\phi(tx+sy)}}.
$$
We apply Pr\'ekopa-Leindler theorem to the functions
$$
u(x)=\exp\PAR{\frac{g(x)}{t}-\phi(x)},\quad v(y)=\exp\PAR{-\phi(y)},\quad w(z)=\exp\PAR{g_s(z)-\phi(z)},
$$
to get
$$
\PAR{\int \exp(g/t)d\mu_\phi}^t\leq\int \exp(g_s)d\mu_\phi.
$$
The derivation of the $L^p$ norm gives the entropy, then using
Taylor formula we get
$$
\PAR{\int \exp(g/t)d\mu_\phi}^t=\int e^g\mu_\phi+s\ent{\mu_\phi}{e^g}+O(s^2).
$$
Then apply Lemma~\ref{lem-technique} to get
\begin{multline*}
\int \exp(g_s)d\mu_\phi=\\
\int e^g\mu_\phi+ s\int
\BRA{z\cdot\nabla g(z)-\phi^*\PAR{\nabla\phi(z)}+
\phi^*\PAR{\nabla\phi(z)-\nabla g(z)}}e^{g(z)} d\mu_\phi(z)+O(s^2).
\end{multline*}
Then when $s$ goes to 0 we get inequality~\eqref{eq-theo}.

Then we can extend the inequality~\eqref{eq-theo} for all function
$g$ smooth enough such that integrals exist. {\proofend\\}

Remark that if $\phi(x)=\NRM{x}^2/2+(n/2)\log(2\pi)$ we obtain the
classical logarithmic Sobolev of Gross for the canonical Gaussian
measure on $\dR^n$.

\subsection{Remarks and examples} \label{sec-ex}

In the  next corollary we give the classical result of
perturbation. Of course we lost the optimal constant given in
inequality~\eqref{eq-theo}. 

If $\Phi$ is a function on $\dR^n$ such that $\int e^{-\Phi} dx<\infty$ we note the probability
 measure $\mu_\Phi$  by 
\begin{equation}
\label{e-phi}
d\mu_\Phi(x)=\frac{e^{-\Phi(x)}}{Z_\Phi}dx,
\end{equation}
where $Z_\Phi=\int {e^{-\Phi(x)}}dx$·

\begin{ecor}
\label{co-pe} Assume that $\phi$ is a $\mathcal C^1$, strictly
convex function on $\dR^n$ such that
$\lim_{\ABS{x}\rightarrow\infty}{\phi(x)}/{\NRM{x}}=\infty$.
Let $\Phi=\phi+U$, where $U$ is a bounded function on $\dR^n$ and
denote by $\mu_\Phi$ the measure defined by~\eqref{e-phi}.

Then for all smooth enough function $g$ on $\dR^n$. we get
\begin{equation}
\label{eq-core} \ent{\mu_\Phi}{e^g}\leq e^{2\text{osc}(U)}\int
\BRA{x\cdot\nabla
g(x)-\phi^*(\nabla\phi(x))+\phi^*\PAR{\nabla\phi(x)-\nabla
g(x)}}e^{g(x)}d\mu_\Phi(x),
\end{equation}
where $\text{osc}(U)=\sup(U)-\inf(U)$.
\end{ecor}

\begin{eproof}
First we observe that
\begin{equation}
\label{eq-ra}
e^{-\text{osc}(U)}\leq\frac{d{\mu}_\Phi}{d\mu_\phi}\leq
e^{\text{osc}(U)}.
\end{equation}
Moreover we have  for all probability measure $\nu$ on $\dR^n$,
$$
\ent{\nu}{e^g}=\inf_{a\geq0}\BRA{\int\PAR{e^g\log\frac{e^g}{a}-e^g+a}d\nu},
$$
using the fact that $\forall x,a>0$, $x \log\frac{x}{a}-x+a\geq0$
we get
$$
e^{-\text{osc}(U)}\ent{{\mu}_\Phi}{e^g}\leq\ent{\mu_\phi}{e^g}\leq
e^{\text{osc}(U)}\ent{{\mu}_\Phi}{e^g}.
$$
Then if $g$ a smooth enough function $g$ on $\dR^n$ we have
\begin{equation*}
\begin{array}{rl}
\disp\ent{{\mu}_\Phi}{e^g}& \disp\leq e^{\text{osc}(U)}\ent{{\mu}_\phi}{e^g}\\
&\disp \leq  e^{\text{osc}(U)} \int \BRA{x\cdot\nabla
g(x)-\phi^*(\nabla\phi(x))+\phi^*\PAR{\nabla\phi(x)-\nabla
g(x)}}e^{g(x)}d\mu_\phi(x).
\end{array}
\end{equation*}
Using the fact that $\phi^*$ is a
convex  function on $\dR^n$ and
$\nabla\phi^*\PAR{\nabla\phi(x)}=x$ we obtain that
$$
\forall x\in\dR^n,\quad {x\cdot\nabla
g(x)-\phi^*(\nabla\phi(x))+\phi^*\PAR{\nabla\phi(x)-\nabla
g(x)}}\geq0.
$$
Then by~\eqref{eq-ra} we get
\begin{equation*}
\ent{{\mu}_\Phi}{e^g} \leq  e^{2\text{osc}(U)} \int \BRA{x\cdot\nabla
g(x)-\phi^*(\nabla\phi(x))+\phi^*\PAR{\nabla\phi(x)-\nabla
g(x)}}e^gd\mu_\Phi.
\end{equation*}
\end{eproof}

\begin{erem}
It is not necessary to give a tensorisation result because we will
obtain exactly the same expression if we compute directly with a
product measure.
\end{erem}

Using Theorem~\ref{the-theo} we  find also the same examples given
in \cite{bobkov-ledoux1} and \cite{bobkov-zeg}.
\begin{ecor}
\label{bl-prop} Let $p\geq 2$ and let $\Phi(x)=\NRM{x}^p/p$ where
$\NRM{\cdot}$ is Euclidean norm in $\dR^n$. Then we get for all
smooth enough function $g$,
\begin{equation}
\label{e-prop} \ent{\mu_\Phi}{e^g}\leq c\int \NRM{\nabla
g}^qe^gd\mu_\Phi,
\end{equation}
where $1/p+1/q=1$ and for some constant $c>0$.
\end{ecor}
\begin{eproof}
Using Theorem~\ref{the-theo}, we just have  to prove that
\begin{equation*}
\forall x\in\dR^n,\quad\forall y\in\dR^n,\quad {x\cdot\nabla
g(x)-\phi^*(\nabla\phi(x))+\phi^*\PAR{\nabla\phi(x)-\nabla
g(x)}}\leq c\NRM{y}^q.
\end{equation*}
Assume that $y\neq 0$ and let note by
$$
\psi(x,y)=\frac{{x\cdot\nabla
g(x)-\phi^*(\nabla\phi(x))+\phi^*\PAR{\nabla\phi(x)-\nabla
g(x)}}}{\NRM{y}^q}.
$$
Then $\psi $ is a bounded function. Indeed
 an easy calculus prove that $\phi^*(x)=\NRM{x}^q/q$. Let take
now $z=x\NRM{x}^{p-2}\NRM{y}$ and $e=y/\NRM{y}$ then we obtain
$$
\psi(x,y)=\bar{\psi}(z,e)=z\cdot
e\NRM{z}^{q-2}-\frac{1}{q}\NRM{z}^q+\frac{1}{q}\NRM{\frac{z}{\NRM{z}}-\frac{e}{\NRM{z}}}^q.
$$
We have $\NRM{e}=1$, then $e$ is bounded. Using Taylor formula we
get $\bar{\psi}(z,e)=O(\NRM{y}^{q-2})$. But $p\geq 2$ implies that
$q\leq2$ and then $\bar{\psi}$ is a bounded function. $\psi$ is
then a bounded, if $c=\sup{\psi}$ we get then
inequality~\eqref{e-prop}.
\end{eproof}

We can remark that Proposition~\ref{bl-prop} is not true when
$p\in]1,2[$. As we can see in \cite{ge-gu-mi}, when $p\in]1,2[$ we
have to change the right hand term of inequality~\eqref{e-prop}
and to add a quadratic term.

\medskip

In  Proposition~2.1 of \cite{bobkov-ledoux1}, Bobkov and Ledoux
prove that Pr\'ekopa-Leindler's theorem implies Brascamp-Lieb
inequality. In our case we prove that Theorem~\ref{the-theo}
implies also some Brascamp-Lieb inequality as we can see in the
next corollary.

\begin{ecor}
Let $\phi$ satisfying conditions of Theorem~\ref{the-theo} and assume that $\phi$ is
$\mathcal C^2$ on $\dR^n$.  Then
for all smooth enough function $g$ we get
\begin{equation*}
\var{\mu_\phi}{g}\leq \int \nabla
g\cdot\text{Hess}(\phi)^{-1}\nabla gd\mu_\phi,
\end{equation*}
where $\text{Hess}(\phi)^{-1}$ denote the inverse of the Hessian
of $\phi$.
\end{ecor}

\begin{eproof}
Assume that $g$ is a ${\mathcal C}^\infty$  function with a
compact support  and let apply inequality~\eqref{eq-theo} with the
function $\ep g$ where $\ep>0$. Using Taylor formula we get
$$
\ent{\mu_\phi}{\exp{\ep f}}=2\ep^2\var{\mu_\phi}{f}+o(\ep^2),
$$
and
\begin{multline*}
\int \BRA{x\cdot\nabla
g(x)-\phi^*(\nabla\phi(x))+\phi^*\PAR{\nabla\phi(x)-\nabla
g(x)}}e^g{(x)}d\mu_\phi(x) =
\\\int\frac{\ep^2}{2}\nabla
g\cdot\text{Hess}(\phi^*)\PAR{\nabla \phi}\nabla
gd\mu_\phi+o(\ep^2).
\end{multline*}
Using the fact that $\nabla\phi^*(\nabla\phi(x))=x$ we get that
$\text{Hess}(\phi^*)\PAR{\nabla\phi}=\text{Hess}(\phi)^{-1}$ and
the corollary is proved.
\end{eproof}

\begin{erem}
Let $\phi$ satisfying properties of Theorem~\ref{the-theo}. Note
$$
L(x,y)=\phi(y)-\phi(x)+(y-x)\nabla\phi(x),
$$
due to the convexity of $\phi$ we get that $L(x,y)\geq0$ for all
$x,y\in\dR^n$.

Let $F$ be a density of probability with respect to the measure
$\mu_\phi$, we defined the following Wasserstein distance with the
cost function equal to $L$ by
$$
 W_L(Fd\mu_\phi,d\mu_\phi)=\inf\BRA{ \int L(x,y)d\pi(x,y)},
$$
where the infimum is taken for all probabilities measures $\pi$ on
$\dR^n\times\dR^n$ with marginal distributions $Fd\mu_\phi$ and
$d\mu_\phi$. Then Bobkov and Ledoux prove again in \cite{bobkov-ledoux1}
the following transportation inequality
\begin{equation}
\label{eq-t} W_L(Fd\mu_\phi,d\mu_\phi)\leq \ent{\mu_\phi}{F}.
\end{equation}

The main theorem of Otto and Villani in \cite{villani} is the
following: Classical logarithmic Sobolev inequality  (when
$\phi(x)=\NRM{x}^2/2+(n/2)\log (2\pi)$) implies the transportation
inequality~\eqref{eq-t} for all function $F$, density of
probability with respect to $\mu_\phi$ (see also \cite{bgl} for an
other proof). By the method developed in~\cite{bgl}, one can
easily extend the property  for $\phi(x)=\NRM{x}^p+Z_\phi$ ($p\geq2$).

In the general case exposed here, we don't know if
inequality~\eqref{eq-theo} imply inequality~\eqref{eq-t}.
\end{erem}

\section{Application to Euclidean logarithmic Sobolev inequality}
\label{sec-e}

\begin{ethm}
\label{theo2} Assume that the function $\phi$ satisfies conditions
of Theorem~\ref{the-theo} then for all $\la>0$ and for all smooth
enough function $g$ on $\dR^n$ such that integrals exits we get
\begin{equation}
\label{e-theo2}
\ent{dx}{e^g}\leq -n\log\PAR{\la e}\int e^gdx+\int\phi^*\PAR{-\la\nabla g}e^gdx.
\end{equation}

Last inequality is optimal in the sense that if $g=-C(x-\bar{x})$
with $\bar{x}\in\dR^n$ and $\la=1$ we get an equality.
\end{ethm}

\begin{eproof}
Using  integration by parts on the second term of~\eqref{eq-theo} we obtain for all $g$ smooth enough
$$
\int {x\cdot\nabla g(x)}e^{g(x)}d\mu_\phi(x)=
\int \PAR{-n+x\cdot\nabla\phi(x)}e^{g(x)}d\mu_\phi(x).
$$
Then using the equality $\phi^*\PAR{\nabla\phi}=x\cdot\nabla\phi(x)-\phi(x)$ we get
for all smooth enough $g$
\begin{equation*}
\ent{\mu_\phi}{e^g}\leq
\int \PAR{-n+\phi+\phi^*\PAR{\nabla\phi-\nabla g}}e^{g}d\mu_\phi,
\end{equation*}
Let now take $g=f+\phi$ to obtain
\begin{equation*}
\ent{dx}{e^f}\leq
\int \PAR{-n+\phi^*\PAR{-\nabla g}}e^{g}dx.
\end{equation*}
Let $\la>0$ and take $f(x)=g(\la x)$ we get then
\begin{equation*}
\ent{dx}{e^g}\leq
-n\log\PAR{\la e}\int e^gdx+\int\phi^*\PAR{-\la\nabla g}e^gdx,
\end{equation*}
which prove~\eqref{e-theo2}.

If now  $g=-C(x-\bar{x})$ with $\bar{x}\in\dR^n$ an easy calculus prove that if $\la=1$ we get an equality.
 \end{eproof}

In the inequality~\eqref{e-theo2}, there exits an optimal
$\la_0>0$. Unfortunately, in the almost case we can't give the
expression of the optimal $\la_0$. It is the unique real
satisfying the following equality
$$
-n\int e^gdx+\la_0\int\nabla g\cdot\nabla(\phi^*)\PAR{-\la_0\nabla
g}e^gdx=0.
$$
But when $C$ is homogeneous, we can give an better expression of the
last theorem. We find inequality called Euclidean logarithmic
Sobolev inequality which is explained  on the next  corollary.

\begin{ecor}
\label{cor1}
Let $C$ a strictly convex function on $\dR^n$ and assume that $C$ is  $q$-homogeneous,
$$\forall \la\geq0\quad \text{and}\quad \forall x\in\dR^n,\quad C(\la x)={\la}^qC(x).$$
Then for all smooth enough function $g$ in $\dR^n$ we get
\begin{equation}
\label{e-ecl} \ent{dx}{e^g}\leq \frac{n}{p}\int
e^gdx\log\PAR{\frac{p}{n e^{p-1}{\mathcal L}^{p/n}}\frac{\int
C^*\PAR{-\nabla g}e^gdx}{\int e^gdx} },
\end{equation}
where ${\mathcal L}=\int e^{-C}dx$ and $1/p+1/q=1$.
\end{ecor}

\begin{eproof}
Let apply Theorem~\ref{theo2} with $ \phi=C+\log{\mathcal L}$. Then $\phi$ satisfies conditions
of Theorem~\ref{theo2} and we get then
\begin{equation*}
\ent{dx}{e^g}\leq -n\log\PAR{\la e {\mathcal L}^{1/n}}\int
e^gdx+\int C^*\PAR{-\la\nabla g}e^gdx.
\end{equation*}
Due to the fact that $C$ is $q$-homogeneous an easy calculus prove
that  $C^*$ is  $p$-homogeneous where $1/p+1/q=1$. An optimization
over $\la>0$ gives inequality~\eqref{e-ecl}.
\end{eproof}

Inequality~\eqref{e-ecl} is called Euclidean logarithmic Sobolev
inequality. This inequality with $p=2$ appears in the work of
Weissler in~\cite{weissler}. It was discussed and extended to this
last version in many articles see
\cite{carlen,ledoux96,beckner,del-dol,gentil03,ghoussoub}.

\begin{erem}
Of course as it is explained in the introduction, calculus used in Corollary~\ref{cor1} prove that
inequality~\eqref{e-ecl} is equivalent to inequality~\eqref{e-theo2}.  Agueh, Ghoussoub and Kang, in~\cite{ghoussoub},
used Monge-Kantorovich theory for mass transport to prove inequalities~\eqref{e-theo2} and~\eqref{e-ecl}. Then it
gives an other way to establish Theorem~\ref{the-theo}.

\medskip

Note also that inequality~\eqref{e-ecl} is optimal, extremal functions is given by
$g(x)=-bC(x-\bar{x})$, with $\bar{x}\in\dR^n$ and $b>0$. But we don't know if it's  only 
extremal functions.
\end{erem}

\section{Application to logarithmic Sobolev inequality for large entropy}
\label{sec3}

In \cite{ge-gu-mi,ge-gu-mi2} is given
 a convex logarithmic Sobolev inequality for measure $\mu_\phi$ between $e^{-\ABS{x}}$ and
$e^{-x^2}$. More precisely let $\Phi$ a function on the real line and assume that $\Phi$ is 
even and satisfies the
following property, there exists $M\geq 0$
 and $0<\e\leq1/2$ such that
\begin{equation*}
\tag{\bf H} \forall{x\geq M},\,\,\,\,(1+\e)\Phi(x)\leq
x\Phi'(x)\leq(2-\e)\Phi(x).
\end{equation*}

Then there exists $A,B,D>0$ such that for all smooth functions $g$
 we have
\begin{equation}
\label{ggl}
 \ent{\mu_\Phi}{e^g}\leq
A\int H_{\Phi}\PAR{g'}e^gd\mu_\Phi,
\end{equation}
where
\begin{equation*}
\label{defh}
H_{\Phi}(x)=\left\{
\begin{array}{rl}
\Phi^*\PAR{Bx} &\text{ if }\ABS{x}\geq D,\\
x^2 &\text{ if }\ABS{x}\leq D,
\end{array}
\right.
\end{equation*}
and $\mu_\Phi$ is defined on~\eqref{e-phi}.

The proof of inequality~\eqref{ggl} is technical and it
divided between two parts: the large and the small entropy. We give in the
next theorem a $n$-dimensional version of this inequality but only
for large entropy.

\begin{ethm}
\label{thmn} Let $\Phi$ be a $\mathcal C^1$, strictly convex and
even function on $\dR^n$, such that $
\lim_{\ABS{x}\rightarrow\infty}{\Phi(x)}/{\NRM{x}}=\infty. $
Assume that $\Phi\geq0$ and $\Phi(0)=0$ {$($}it implies that $0$ is the unique
minimum of $\Phi${$)$}. 

Assume that
\begin{equation}
\label{eqm}
\lim_{\al\rightarrow0,\,\al\in[0,1]}
\sup_{x\in\dR^n}\BRA{(1-\al)\frac{\Phi^*\PAR{\frac{x}{1-\al}}}{\Phi^*(x)}}=1,
\end{equation}
assume also that there exists $A>0$ such that
\begin{equation}
\label{eqh}
\forall x\in\dR^n,\quad x\cdot \nabla \Phi(x)\leq (A+1)\Phi(x).
\end{equation}

Then there exists $C_1,C_2\geq0$ such that for all smooth enough function $g$ such that
$\int e^gd\mu_\Phi=1$ and $\ent{\mu_\Phi}{e^g}\geq1$ we get
\begin{equation}
\label{eqd}
\ent{\mu_\Phi}{e^g}\leq C_1\int\Phi^*\PAR{C_2\nabla g}e^gd\mu_\Phi.
\end{equation}
\end{ethm}

\begin{eproof}
Let apply Theorem~\ref{the-theo} with $\phi=\Phi+\log Z_\Phi$ we get then
\begin{equation*}
\ent{\mu_\Phi}{e^g}\leq \int \BRA{x\cdot\nabla
g(x)-\Phi^*(\nabla\Phi(x))+\Phi^*\PAR{\nabla\Phi(x)-\nabla
g(x)}}e^gd\mu_\Phi.
\end{equation*}

Let $\al\in[0,1[$, $\Phi^*$ is convex then
\begin{equation}
\label{eqc}
\Phi^*\PAR{\nabla\Phi(x)-\nabla g(x)}\leq(1-\al)\Phi^*\PAR{\frac{\nabla\Phi(x)}{1-\al}}+
\al\Phi^*\PAR{\frac{-\nabla g(x)}{\al}},
\end{equation}
recall that $\Phi^*$ is also a even function.
Young's inequality implies that
\begin{equation}
\label{eqy}
x\cdot\frac{\nabla g(x)}{\al}\leq\Phi(x)+\Phi^*\PAR{\frac{\nabla g(x)}{\al}}.
\end{equation}
Using~\eqref{eqc} and ~\eqref{eqy} we get
\begin{multline*}
\ent{\mu_\Phi}{e^g}\leq 2\al\int \Phi^*\PAR{\frac{\nabla g}{\al}}e^gd\mu_\Phi+
\al \int \Phi\PAR{x}e^gd\mu_\Phi+\\
\int \PAR{(1-\al)\Phi^*\PAR{\frac{\nabla\Phi(x)}{1-\al}}-\Phi^*\PAR{{\nabla\Phi(x)}}}e^gd\mu_\Phi.
\end{multline*}
We have $\Phi^*(\nabla \Phi(x))=x\cdot\nabla \Phi(x)-\Phi(x)$, then
inequality~\eqref{eqh} implies that $\Phi^*\PAR{\nabla \Phi(x)}\leq A\Phi(x)$. Due to the fact that
$\Phi(0)=0$ we have $\Phi^*\geq 0$ we get
$$
\ent{\mu_\Phi}{e^g}\leq \al\int \Phi^*\PAR{\frac{\nabla g}{\al}}e^gd\mu_\Phi+
 \al\int \Phi^*\PAR{\frac{\nabla g}{\al}}e^gd\mu_\Phi+
(\al+A\ABS{\psi(\al)-1})\int \Phi e^gd\mu_\Phi,
$$
where
\begin{equation}
\label{eqp}
\psi(\al)=\sup_{x\in\dR^n}\BRA{(1-\al)\frac{\Phi^*\PAR{\frac{x}{1-\al}}}{\Phi^*(x)}}.
\end{equation}
Let $\la>0$ then due to the fact that $\int e^g d\mu_\Phi=1$ we get
$$
\int \Phi e^gd\mu_\Phi\leq \la\PAR{\ent{\mu_\Phi}{e^g}+\log\int e^{\Phi /\la}d\mu_\Phi}.
$$
We have $\disp\lim_{\la\rightarrow\infty}\log\int e^{\Phi /\la}d\mu_\Phi=0$, then
let now choose $\la$ large enough such that  $\log\int e^{\Phi /\la}d\mu_\Phi\leq 1$.
Using the property~\eqref{eqm}, take $\al$ such that
$(\al+A\ABS{\psi(\al)-1})\la\leq1/4$ we obtain
$$
\ent{\mu_\Phi}{e^g}\leq 2\al\int \Phi^*\PAR{\frac{\nabla g}{\al}}e^gd\mu_\Phi+
\frac{1}{4}\PAR{\ent{\mu_\Phi}{e^g}+1}.
$$
Then using $\ent{\mu_\Phi}{e^g}\geq1$ we obtain
$$
\ent{\mu_\Phi}{e^g}\leq{4\al}\int \Phi^*\PAR{\frac{\nabla g}{\al}}e^gd\mu_\Phi.
$$
\end{eproof}

We need a lemma to  give non-trivial examples. This  lemma
explains how property~\eqref{eqm} is a infinity property.

\begin{elem}
\label{lm1}
Let $\Phi_1$ and $\Phi_2$ be two strictly convex and even functions such that 
$\Phi_1,\Phi_2\geq0$, $\Phi_1(0)=\Phi_2(0)=0$
and
$
\lim_{\ABS{x}\rightarrow\infty}{\Phi_1(x)}/{\NRM{x}}=
\lim_{\ABS{x}\rightarrow\infty}{\Phi_2(x)}/{\NRM{x}}=\infty.
$
Assume also that
$\disp\Phi_1(x)\stackrel{{\pm\infty}}{\sim}\Phi_2(x)$. 

\bigskip

If $\Phi_2$
satisfies the property~\eqref{eqm} then $\Phi_1$ satisfies also
the same property.
\end{elem}

\begin{eproof}
First we prove that  $\disp\Phi_1^*(x)\stackrel{{\pm\infty}}{\sim}\Phi_2^*(x)$.
Let $\ep>0$, then there exists $A>0$ such that
$$
\forall y\in\dR^n,\quad \NRM{y}\geq A,\quad (1-\ep)\Phi_2(y)\leq \Phi_1(y)\leq(1+\ep)\Phi_2(y),
$$
then
$$
\forall x\in\dR^n,\quad \sup_{\NRM{y}\geq A}\BRA{x\cdot y-(1+\ep)\Phi_2(y)}\leq
\sup_{\NRM{y}\geq A}\BRA{x\cdot y-\Phi_1(y)}\leq \sup_{\NRM{y}\geq A}\BRA{x\cdot y-(1-\ep)\Phi_2(y)}.
$$
$\Phi_1$ and $\Phi_2$ are strictly convex then there exists $B>0$
such that 
$$
\forall x\in\dR^n,\quad\NRM{x}\geq B,\quad \Phi_1^*(x)=\sup_{\NRM{y}\geq A}\BRA{x\cdot y-\Phi_1(y)},
$$
and the same for $\Phi_2$, then  
$$
\forall x\in\dR^n,\, \NRM{x}\geq B,\quad (1+\ep)\Phi_2^*\PAR{\frac{x}{1+\ep}}\leq \Phi_1^*(x)\leq(1-\ep)\Phi_2^*\PAR{\frac{x}{1-\ep}}.
$$

Using now property~\eqref{eqm} for $\Phi_2$ we get
$$
\forall x\in\dR^n,\quad \Phi^*_2\PAR{x}\leq {\psi\PAR{\frac{\ep}{1+\ep}}}\PAR{1+\ep}\Phi_2^*\PAR{\frac{x}{1+\ep}}
\text{ and }
 \Phi^*_2\PAR{\frac{x}{1-\ep}}\leq\frac{\psi(\ep)}{1-\ep}\Phi_2^*\PAR{x},
$$
where $\psi$ is defined on~\eqref{eqp}. We get then
$$
\forall  x\in\dR^n,\, \NRM{x}\geq B,\quad
\psi\PAR{\frac{\ep}{1+\ep}}^{-1}\Phi_2^*\PAR{{x}}\leq
\Phi_1^*(x)\leq\psi(\ep)\Phi_2^*\PAR{x}.
$$
The function $\Phi_2$ satisfies~\eqref{eqm} then
$\lim_{\al\rightarrow0}\psi(\al)=1$ then
$\disp\Phi_1^*(x)\stackrel{{\pm\infty}}{\sim}\Phi_2^*(x)$.

The end of the proof is elementary, we  just have to remark that using a compact argument we get
$$
\forall A>0,\quad \lim_{\al\rightarrow0,\,\al\in[0,1]}
\sup_{\NRM{x}\leq A}\BRA{(1-\al)\frac{\Phi_1^*\PAR{\frac{x}{1-\al}}}{\Phi_1^*(x)}}=1.
$$
Then, when $\NRM{x}$ is large  $\Phi_1^*$ is equivalent to $\Phi_2^*$.
\end{eproof}

\begin{exe}
\begin{itemize}
\item Let $\Phi$ be a $\mathcal C^1$, strictly convex function on $\dR^1$.
Assume that $\Phi\geq0$ and $\Phi(0)=0$. Assume that
$$
\forall x\in\dR,\,\ABS{x}\geq2,\quad
\Phi(x)=x^a\log^b x
$$ with
$a>1$ and $b\in\dR$. Then $\Phi$ satisfies property~\eqref{eqm}.
Remark that if $a\in]1,2[$ and $b=0$ then the measure $\mu_{\Phi}$
doesn't satisfies~\eqref{e-prop} for small entropy.
\item  Here is 
 now an example of  measure on $\dR^n$ with interactions. Let $h$ be a $\mathcal C^1$,
strictly convex function on $\dR^1$.
Assume that $h\geq0$, $h(0)=0$ and that
$h$ satisfies assumptions~\eqref{eqm} and~\eqref{eqh}. Assume also that
\begin{equation}
\label{eqo}
\lim_{\ABS{x}\rightarrow\infty}\frac{h(x)}{x^2}=+\infty.
\end{equation}
Note
$$
\Phi(x)=\sum_{i=1}^n (x_ix_{i+1}+h(x_i)),
$$
where $x=(x_1,\cdots,x_n)$ and $x_{n+1}=x_1$. Then it's easy to prove that $\Phi$ is convex, even with $\Phi(0)=0$ and
satisfies inequality~\eqref{eqh}. Then using~\eqref{eqo} we get that
$$\Phi(x)\stackrel{{\pm\infty}}{\sim}\sum_{i=1}^n h(x_i).$$ 

By Lemma~\ref{lm1} we prove that
$\Phi$ satisfy~\eqref{eqm}.

This example in interesting because it gives an measure on $\dR^n$ which is not a product measure on $\dR^n$ and satisfies
inequality~\eqref{eqd} for large entropy.
\end{itemize}
\end{exe}

\newcommand{\etalchar}[1]{$^{#1}$}


\begin{thebibliography}{GGM05b}

\bibitem[ABC{\etalchar{+}}00]{logsob}
C.~An{\'e}, S.~Blach{\`e}re, D.~Chafa{\"\i}, P.~Foug{\`e}res, I.~Gentil,
  F.~Malrieu, C.~Roberto, and G.~Scheffer.
\newblock {\em Sur les in{\'e}galit{\'e}s de {S}obolev logarithmiques},
  volume~10 of {\em Panoramas et Synth{\`e}ses}.
\newblock Soci{\'e}t{\'e} {M}ath{\'e}matique de {F}rance, Paris, 2000.

\bibitem[AGK04]{ghoussoub}
M.~Agueh, N.~Ghoussoub, and X.~Kang.
\newblock {Geometric inequalities via a general comparison principle for
  interacting gases.}
\newblock {\em Geom. Funct. Anal.}, 14(1):215--244, 2004.

\bibitem[Bec99]{beckner}
W.~Beckner.
\newblock {Geometric asymptotics and the logarithmic Sobolev inequality.}
\newblock {\em Forum Math.}, 11(1):105--137, 1999.

\bibitem[BGL01]{bgl}
S.~Bobkov, I.~Gentil, and M.~Ledoux.
\newblock Hypercontractivity of {H}amilton-{J}acobi equations.
\newblock {\em J. Math. Pu. Appli.}, 80(7):669--696, 2001.

\bibitem[BL00]{bobkov-ledoux1}
S.~G. Bobkov and M.~Ledoux.
\newblock From {B}runn-{M}inkowski to {B}rascamp-{L}ieb and to logarithmic
  {S}obolev inequalities.
\newblock {\em Geom. Funct. Anal.}, 10(5):1028--1052, 2000.

\bibitem[BZ05]{bobkov-zeg}
S.~Bobkov and B.~Zegarlinski.
\newblock Entropy bounds and isoperimetry.
\newblock To appear in {M}emoirs {A}{M}{S}, 2005.

\bibitem[Car91]{carlen}
E.~A. Carlen.
\newblock {Superadditivity of Fisher's information and logarithmic Sobolev
  inequalities.}
\newblock {\em J. Funct. Anal.}, 101(1):194--211, 1991.

\bibitem[DPD03]{del-dol}
M.~Del~Pino and J.~Dolbeault.
\newblock {The optimal Euclidean $L^{p}$-Sobolev logarithmic inequality.}
\newblock {\em J. Funct. Anal.}, 197(1):151--161, 2003.

\bibitem[Gen03]{gentil03}
I.~Gentil.
\newblock The general optimal {$L\sp p$}-{E}uclidean logarithmic {S}obolev
  inequality by {H}amilton-{J}acobi equations.
\newblock {\em J. Funct. Anal.}, 202(2):591--599, 2003.

\bibitem[GGM05a]{ge-gu-mi2}
I.~Gentil, A.~Guillin, and L.~Miclo.
\newblock Logarithmic sobolev inequalities in curvature null.
\newblock In preparation, 2005.

\bibitem[GGM05b]{ge-gu-mi}
I.~Gentil, A.~Guillin, and L.~Miclo.
\newblock Modified logarithmic sobolev inequalities and transportation
  inequalities.
\newblock To appear in {P}robab. {T}heory {R}elated {F}ields, 2005.

\bibitem[Gro75]{gross}
L.~Gross.
\newblock Logarithmic {S}obolev inequalities.
\newblock {\em Amer. J. Math.}, 97(4):1061--1083, 1975.

\bibitem[Gup80]{gupta}
S.~D. Gupta.
\newblock {Brunn-Minkowski inequality and its aftermath.}
\newblock {\em J. Multivariate Anal.}, 10:296--318, 1980.

\bibitem[Led96]{ledoux96}
M.~Ledoux.
\newblock {Isoperimetry and Gaussian analysis.}
\newblock In {\em {Dobrushin, R. (ed.) et al., Lectures on probability theory
  and statistics. Ecole d'\'et\'e de probabilit\'es de Saint-Flour XXIV -- 1994.
  Berlin: Springer. Lect. Notes Math. 1648, 165-294 }}. 1996.

\bibitem[Mau04]{maurey}
B.~Maurey.
\newblock In{\'e}galit{\'e} de {B}runn-{M}inkowski-{L}usternik, et autres
  in{\'e}galit{\'e}s g{\'e}om{\'e}triques et fonctionnelles.
\newblock {\em S{\'e}minaire Bourbaki}, 928, 2003/04.

\bibitem[OV00]{villani}
F.~Otto and C.~Villani.
\newblock Generalization of an inequality by {T}alagrand, and links with the
  logarithmic {S}obolev inequality.
\newblock {\em J. Funct. Anal.}, 173(2):361--400, 2000.

\bibitem[Wei78]{weissler}
F.~B. Weissler.
\newblock {Logarithmic Sobolev inequalities for the heat-diffusion semigroup.}
\newblock {\em Trans. Am. Math. Soc.}, 237:255--269, 1978.

\end{thebibliography}
\end{document}